\documentclass[a4paper,12pt]{article}

\usepackage{amsmath, amssymb, amsthm}

\usepackage[dvips,final]{graphicx}

\usepackage{color}

\usepackage{titlesec}
\titleformat{\section}{\bfseries}{\thesection}{1em}{}
\titleformat{\subsection}{\itshape}{\thesubsection}{1em}{}

\numberwithin{equation}{section}

\def\ve{\varepsilon}

\def\dd{\,\mathrm{d}}
\def\dive{\mathrm{\,div\,}}
\def\for{\mathrm{for\ }}
\def\supess{\mathop{\rm sup\,ess\,}}
\def\on{^{(n)}}

\def\real{\mathbb{R}}
\def\nat{\mathbb{N}}

\def\io{\int_{\Omega}}
\def\ipo{\int_{\partial\Omega}}

\def\be{\begin{equation}\label}
	\def\ee{\end{equation}}

\def\ber{\begin{eqnarray}}
	\def\eer{\end{eqnarray}}

\def\bers{\begin{eqnarray*}}
	\def\eers{\end{eqnarray*}}

\def\bpf{\begin{pf}}
	\def\epf{\end{pf}}

\newtheorem{theorem}{Theorem}[section]
\newtheorem{lemma}[theorem]{Lemma}
\newtheorem{corollary}[theorem]{Corollary}
\newtheorem{hypothesis}[theorem]{Hypothesis}

\newtheorem{remark}[theorem]{Remark}

\begin{document}
	

\title{A model for lime consolidation of porous solids\thanks{This work was supported by the GA\v CR Grant No.~20-14736S; the European Regional Development Fund Project No. CZ.02.1.01/0.0/0.0/16{\_}019/0000778; and the Austrian Science Fund (FWF) Project V662.
}}

\author{Bettina Detmann
\thanks{University of Duisburg-Essen,
Faculty of Engineering, Department of Civil Engineering,
D-45117 Essen, Germany, E-mail: \tt bettina.detmann@uni-due.de.}
\and Chiara Gavioli
\thanks{Institute of Analysis and Scientific Computing, TU Wien, Wiedner Hauptstra\ss e 8-10, A-1040 Vienna (Austria), E-mail: {\tt chiara.gavioli@tuwien.ac.at}.}
\and Pavel Krej\v c\'{\i}
\thanks{Faculty of Civil Engineering, Czech Technical University, Th\'akurova 7, CZ-16629 Praha 6, Czech Republic, E-mail: {\tt Pavel.Krejci@cvut.cz}.}
\and Jan Lama\v c
\thanks{Faculty of Civil Engineering, Czech Technical University, Th\'akurova 7, CZ-16629 Praha 6, Czech Republic, E-mail: {\tt Jan.Lamac@cvut.cz}.}
\and Yuliya Namlyeyeva
\thanks{Faculty of Civil Engineering, Czech Technical University, Th\'akurova 7, CZ-16629 Praha 6, Czech Republic, E-mail: {\tt yuliya.namlyeyeva@fsv.cvut.cz}.}
}

\date{}

\maketitle

\begin{abstract}
We propose a mathematical model describing the process of filling the pores of a building material with lime water solution with the goal to 
improve the consistency of the porous solid. Chemical reactions produce calcium carbonate which glues the solid particles together at some distance from the boundary and strengthens the whole structure. The model consists of a 3D convection-diffusion system with a nonlinear boundary condition for the liquid and for calcium hydroxide, coupled with the mass balance 
equations for the chemical reaction. The main result consists in proving that the system has a solution for each initial data from a physically relevant class. A 1D numerical test shows a qualitative agreement with experimental observations.

\bigskip


\medskip


\end{abstract}

\section*{Introduction}

In this paper we investigate mathematically a process which is used by the building industry in order to protect and conserve cultural goods and other structure works. Such structures which are subject to weathering can 
be strengthened by filling the pores by a water-lime-mixture. The mixture 
penetrates into the pore structure of the stone and the calcium hydroxide 
reacts with carbon dioxide and builds calcium carbonate and water. A solid layer is built in the pore space which strengthens the material. The main problem is, however, that the consolidated layer is in practice rather thin and is located too close to the active boundary.

In order to avoid possible ambiguity, it is necessary to explain that the term `consolidation' is to be interpreted here as a process of formation of calcium carbonate which has the property of binding the particles together. It might also be called `cementation' or `compaction'.

In the literature there are a lot of works dealing with those problems. Different practical strategies for the wetting and drying regime which lead to a more uniform distribution of the consolidant are compared in \cite{sli}. A mathematical model is proposed in \cite{ssv1993, ssv1995}, and \cite{cetal}, where the authors derive governing equations for moisture, heat, and air flow through concrete. A numerical procedure based on the finite element method is developed there to solve the set of equations and to investigate the influence of relative humidity and temperature. It is shown  that the amount of calcium carbonate formed in a unit of time depends on the degree of carbonation,~i.\,e., the availability of calcium hydroxide, the temperature, the carbon dioxide concentration and the relative humidity in the pore structure of the concrete. An extension of the aforementioned papers by studying the hygro-thermal behavior of concrete in the special situation of high temperatures can be found in \cite{getal}.

In the present case chemical reactions take place. Various approaches exist which describe such processes by models which stem from different backgrounds (e.\,g., from mixture theory or empirical models). An overview on 
the development of theories especially for porous media including chemical reactions is given in \cite{bd}.
The interactions between the constituents of a porous medium are not necessarily of chemical nature which would lead to a chemical transformation of one set of chemical substances to another. Simpler is the mass exchange between the constituents by physical processes like adsorption. Adsorption-diffusion processes have been studied by B. Albers (the former 
name of B. Detmann) e.\,g.~in \cite{ads}. Other works on sorption in porous solids including molecular condensation are \cite{bazant1} or \cite{bazant2}. In these works the diffusivities of water and carbon dioxide are assumed to be strongly dependent on pore humidity, temperature and also on the degree of hydration of concrete. The authors realized that the porosity becomes non-uniform in time. This is an observation which is interesting also in the present case because the structure of the channels clearly changes with the progress of the reaction. A survey of consolidation techniques for historical materials is published in \cite{pd}. The influence of the particle size on the efficiency of the consolidation process is 
investigated in \cite{svmlma}. Experimental determination of the penetration depth is the subject of \cite{blhvs1,blhvs2,blhvs3}. Different variants of the consolidants is studied in \cite{drsli,msm,poam,ss}, and an experimental work on mechanical interaction between the consolidant and the matrix material is carried out in \cite{mpe}.

A further work dealing with chemical reactions and diffusion in concrete based on the mixture theory for fluids introduced by Truesdell and coworkers is by A.\,J.~Vromans et al.~\cite{vetal}. The model describes the corrosion of concrete with sulfuric acid which means a transformation of slaked lime and sulfuric acid into gypsum releasing water. It is a similar reaction we are looking at. A similar topic is dealt with in \cite{VanBalen}, where it is shown how the carbonation process in lime mortar is influenced by the diffusion of carbon dioxide into the mortar pore system by 
the kinetics of the lime carbonation reaction and by the drying and wetting process in the mortar.

Experimental results of $\mathrm{CaCO_{3}}$ precipitation kinetics can be 
found in \cite{Roques}. The porosity changes during the reaction. This was studied by Houst and Wittmann, who also investigated the influence of the 
water content on the diffusivity of $\mathrm{CO_{2}}$ and $\mathrm{O_{2}}$ through hydrated cement paste \cite{Houst}. An investigation of the physico-chemical characteristics of ancient mortars with comparison to a reaction-diffusion model by Zouridakis et al.~is presented in \cite{zour}. A 
slightly different reaction involving also sulfur is mathematically studied in \cite{bohm} by B\"ohm et al. There the corrosion in a sewer pipe is 
modeled as a moving-boundary system. A strategy for predicting the penetration of carbonation reaction fronts in concrete was proposed by Muntean et al. in \cite{munt}. A simple 1D mathematical model for the treatment of sandstone neglecting the effects of chemical reactions is proposed in \cite{cnns} and further refined in \cite{bfngrt}.

We model the consolidation process as a convection-diffusion system coupled with chemical reaction in a 3D porous solid. The physical observation that only water can be evacuated from the porous body, while lime remains 
inside, requires a nonstandard boundary condition on the active part of the boundary. We choose a simple one-sided condition for the lime exchange 
between the interior and the exterior. The main result of the paper consists in proving rigorously that the resulting initial-boundary value problem for the PDE system in 3D has a solution satisfying natural physical constraints, including the boundedness of the concentrations proved by means 
of time discrete Moser iterations. We also show the result of numerical simulation in a simplified 1D situation.

The structure of the present paper is the following. In Section \ref{mod}, we explain the modeling hypotheses and derive the corresponding system of balance equations with nonlinear boundary conditions. In Section \ref{math}, we give a rigorous formulation of the initial-boundary value problem, 
specify the mathematical hypotheses, and state the main result in Theorem \ref{t1}. The solution is constructed by a time-discretization scheme proposed in Section \ref{time}. The estimates independent of the time 
step size derived for this time-discrete system constitute the  substantial step in the proof of Theorem \ref{t1}, which is obtained in Section \ref{limi} by passing to the limit as the time step tends to zero. A numerical test for a reduced 1D system is carried out in Section \ref{nume} to illustrate a qualitative agreement of the mathematical result with experimental observations.

	
	\section{The model}\label{mod}
	
	We imagine a porous medium (sandstone, for example) the structure of which is to be strengthened by letting calcium hydroxide particles driven by 
water flow penetrate into the pores. In contact with the air present in the pores, the calcium hydroxide reacts with the carbon dioxide contained in the air and produces a precipitate (calcium carbonate) which is not water-soluble, remains in the pores, and glues the sandstone particles together. Unlike, e.\,g., in \cite{hkk,schK}, we do not consider the porosity as one of the state variables. The porosity evolution law is replaced with the assumption that the permeability decreases as a result of the calcium carbonate deposit in the pores. The chemical reaction is assumed to be irreversible and we write it 
as $Ca(OH)_2 + CO_2 \to CaCO_3 + H_2O$.\medskip
	
	Notation:
	
	\begin{itemize}
		\item[] $\dot c^W$ ... mass source rate of $H_2 O$ produced by the chemical reaction
		\item[] $\dot c^H$ ... mass source rate of $Ca(OH)_2$ produced by the chemical reaction
		\item[] $\dot c^P$ ... mass source rate of $Ca C O_3$ produced by the chemical reaction
		\item[] $\dot c^G$ ... mass source rate of $C O_2$ produced by the chemical reaction
		\item[] $m^W$ ... molar mass of $H_2 O$
		\item[] $m^H$ ... molar mass of $Ca(OH)_2$
		\item[] $m^P$ ... molar mass of $Ca C O_3$
		\item[] $m^G$ ... molar mass of $C O_2$
		\item[] $\rho^W$ ... mass density of $H_2 O$
		\item[] $\rho^H$ ... mass density of $Ca(OH)_2$
		\item[] $p$ ... capillary pressure
		\item[] $s$ ... water volume saturation
		\item[] ${h}$ ... relative concentration of $Ca(OH)_2$
		\item[] $p^{\partial\Omega}$ ... outer pressure
		\item[] $h^{\partial\Omega}$ ... outer concentration of $Ca(OH)_2$
		\item[] $v$ ... transport velocity vector
		\item[] $k(c^P)$ ... permeability of the porous solid
		\item[] $q$ ... liquid mass flux
		\item[] $q^H$ ... mass flux of $Ca(OH)_2$
		\item[] $\gamma$ ... speed of the chemical reaction
		\item[] $\kappa$ ... diffusivity of $Ca(OH)_2$
		\item[] $n$ ... unit outward normal vector
		\item[] $\sigma(x)$ ... transport velocity interaction kernel
		\item[] $\alpha(x)$ ... boundary permeability for water
		\item[] $\beta(x)$ ... boundary permeability for the inflow of $Ca(OH)_2$
	\end{itemize}
	
	Mass balance of the chemical reaction:
	$$
	\frac{\dot c^P}{m^P} = \frac{\dot c^W}{m^W} = -\frac{\dot c^H}{m^H} = -\frac{\dot c^G}{m^G}.
	$$
	Water mass balance in an arbitrary subdomain $V$ of the porous body:
	$$
	\frac{\dd}{\dd t}\int_V \rho^W s \dd x + \int_{\partial V} q\cdot n \dd S = \int_V \dot c^W \dd x.
	$$
	Calcium hydroxide mass balance in an arbitrary subdomain $V$ of the porous body:
	$$
	\frac{\dd}{\dd t}\int_V \rho^H {h} \dd x + \int_{\partial V} q^H\cdot n \dd S = \int_V \dot c^H \dd x.
	$$
	Water mass balance in differential form:
	$$
	\rho^W\dot s + \dive q = \dot c^W.
	$$
	Calcium hydroxide mass balance in differential form:
	$$
	\rho^H\dot{h} + \dive q^H  = \dot c^H.
	$$
	The water mass flux is assumed to obey the Darcy law:
	$$
	q = -k(c^P) \nabla p,
	$$
	with permeability coefficient $k(c^P)$ which is assumed to decrease as the amount $c^P$ of $Ca C O_3$ given by the formula
	$$
	c^P(x,t) = \int_0^t \dot c^P(x,t')\dd t'
	$$
	increases and fills the pores.
	
	The flux of $Ca(OH)_2$ consists of transport and diffusion terms:
	$$
	q^H = \rho^H {h} v - \kappa s \nabla {h}.
	$$
	The mobility coefficient $\kappa s$ in the diffusion term is assumed to be proportional to $s$: If there is no water in the pores, no diffusion takes place.
	
	We assume that the transport of $Ca(OH)_2$ at the point $x \in \Omega$ is driven by the water flux in a small neighborhood of $x$. In mathematical terms, we assume that there exists a nonnegative function $\sigma$ with 
support in a small neighborhood of the origin such that the transport velocity $v$ can be defined as
	$$
	v(x,t) = \frac{1}{\rho^H}\io \sigma(x-y) q(y,t) \dd y.
	$$
The main reason for this assumption is a mathematical one. The strong nonlinear coupling between $s$ and $h$ makes it difficult to control the bounds for the unknowns in the approximation scheme. We believe that such a regularization of the transport velocity makes physically sense as well.

	The wetting-dewetting curve is described by an increasing function $f$:
	$$p=f(s).$$
We focus on modeling the chemical reactions. Capillary hysteresis, deformations of the solid matrix, and thermal effects are therefore neglected here. We plan to include them following the ideas of \cite{ak} in a subsequent study.

	The dynamics of the chemical reaction is modeled according to the so-called \textit{law of mass action}, which states that the rate of the chemical reaction is directly proportional to the product of the concentrations of the reactants. We assume it in the form
	\be{reac}
	\dot c^P = \gamma m^P {h} s(1-s).
	\ee
	Its meaning is that no reaction can take place if either no $Ca(OH)_2$ is available (that is, $h=0$), or no water is available (that is, $s=0$), or no $C O_2$ is available (that is, $s=1$), according to the hypothesis that the chemical reaction takes dominantly place on the contact between water and air. In order to reduce the complexity of the problem, we assume directly that the available quantity of $C O_2$  is proportional to the air content.
	
	On the boundary $\partial\Omega$ we prescribe the normal fluxes. For the 
normal component of $q$, we assume that it is proportional to the difference between the pressures $p$ inside and $p^{\partial\Omega}$ outside the 
body. For the flux of $Ca(OH)_2$, we assume that it can point only inward 
proportionally to the difference of concentrations and to $s$, and no outward flux is possible. Inward flux takes place only if the outer concentration $h^{\partial\Omega}$ is bigger than the inner concentration $h$:
	$$
	\left.
	\begin{array}{rcl}	
	q\cdot n &=& \alpha(x) (p-p^{\partial\Omega})\\[1.5 mm]
	q^H \cdot n &=& -\beta(x) s (h^{\partial\Omega}- h)^+
	\end{array}
	\right\rbrace \text{on }\partial\Omega.$$

	
	\section{Mathematical problem}\label{math}
	
	Let $\Omega \subset \real^3$ be a bounded Lipschitzian domain. We consider the Hilbert triplet $V \subset H \equiv H' \subset V'$ with compact embeddings and with $H = L^2(\Omega)$, $V= W^{1,2}(\Omega)$. For two unknown functions $s(x,t),h(x,t)$ defined for $(x,t) \in\Omega\times (0,T)$ 
the resulting PDE system reads
	\begin{align}\nonumber
		&\io (\rho^W s_t\phi(x) + k(c^P) f'(s) \nabla s\cdot\nabla \phi(x))\dd x + \ipo \alpha(x)(f(s) - f(s^{\partial\Omega}))\phi(x) \dd S(x)\\ \label{e01}
		&\qquad = \gamma m^W \io{h} s(1-s)\phi(x)\dd x,\\ \nonumber
		&\io (\rho^H h_t\psi(x) + (\kappa s \nabla {h} - \rho^H hv)\cdot\nabla\psi(x))\dd x - \ipo \beta(x)s (h^{\partial\Omega} - h)^+\psi(x)\dd S(x)\\ 
\label{e02}
		&\qquad = -\gamma m^H \io {h} s(1-s) \psi(x)\dd x.
	\end{align}
	for all test functions $\phi, \psi \in V$, where $s^{\partial\Omega} := 
f^{-1}(p^{\partial\Omega})$, and with initial conditions
	\be{ini}
	s(x,0) = s^0(x), \quad h(x,0) = h^0(x) \quad \for x \in \Omega.
	\ee
	
	\begin{hypothesis}\label{h1}
		The data have the properties
		\begin{itemize}
			\item[{\rm (i)}] $s^{\partial\Omega} \in L^\infty(\partial\Omega\times 
(0,T))$, $s^0 \in L^\infty(\Omega)\cap W^{1,2}(\Omega)$ are given such that $s^{\partial\Omega}_t \in L^2(\partial\Omega\times (0,T))$, $0<s^\flat\le s^{\partial\Omega}(x,t) \le 1$ for a.\,e. $(x,t) \in \partial\Omega {\times} (0,T)$, $0<s^\flat \le s^0(x) \le 1$ for a.\,e. $x \in \Omega$;
			\item[{\rm (ii)}] $h^{\partial\Omega} \in L^\infty(\partial\Omega\times (0,T))$, $h^0 \in L^\infty(\Omega)$ are given such that  $0\le h^{\partial\Omega}(x,t) \le h^\sharp$ for a.\,e. $(x,t) \in \partial\Omega \times (0,T)$, $0 \le h^0(x) \le h^\sharp$ for some $h^\sharp > 0$ and for a.\,e. $x \in \Omega$;
			\item[{\rm (iii)}] $f:[0,1] \to \real$ is continuously differentiable, 
$0 < f^\flat \le f'(s) \le f^\sharp$ for $0\le s \le 1$;
			\item[{\rm (iv)}] $k$ is continuously differentiable and nonincreasing, $0 < k^\flat \le k(r) \le k^\sharp$ for $r\ge 0$;
			\item[{\rm (v)}] $\sigma: \real^3 \to [0,\infty)$ is continuous with compact support, $\int_{\real^3}\sigma(x)\dd x = 1$;
			\item[{\rm (vi)}] $\alpha, \beta \in L^\infty(\partial\Omega)$, $\alpha(x)\ge 0$, $\beta(x) \ge 0$ on $\partial\Omega$, $\ipo \alpha(x)\dd S(x) > 0$, $\ipo \beta(x)\dd S(x) > 0$.
		\end{itemize}
	\end{hypothesis}

The meaning of Hypothesis \ref{h1}\,(vi) is that the boundary $\partial\Omega$ is inhomogeneous, with different permeabilities at different parts of the boundary. The transport of water (supply of $Ca(OH)_2$) through the boundary takes place only on parts of $\partial\Omega$ where $\alpha >0$ ($\beta>0$, respectively).

The remaining sections are devoted to the proof of the following result.
	
	\begin{theorem}\label{t1}
		Let Hypothesis \ref{h1} hold. Then system \eqref{e01}--\eqref{e02} with initial conditions \eqref{ini} admits a solution $(s,h)$ such that 
$s_t \in L^2(\Omega\times (0,T))$, $\nabla s \in L^\infty(0,T;L^2(\Omega;\real^3))$, $\nabla h \in L^2(\Omega\times(0,T);\real^3)$, $h_t \in L^2(0,T; W^{-1,2}(\Omega))$, $h \in L^\infty(\Omega\times (0,T))$, $s(x,t) \in [0,1]$ a.\,e., $h(x,t) \ge 0$ a.\,e.
	\end{theorem}

We omit the positive physical constants which are not relevant for the analysis. The strategy of the proof is based on choosing a cut-off parameter $R>0$, replacing $h$ in the nonlinear terms with $Q_R(h) = \min\{h^+, R\}$, $1-s$ with $Q_R(1-s)$, and $v$ in \eqref{e02} with $v^R := (Q_R(|v|)/|v|)v$. We also extend the values of the function $f$ outside the interval $[0,1]$ by introducing the function $\tilde f$ by the formula
	$$
	\tilde f(s) = \left\{
	\begin{array}{ll}
		f(0) + f'(0) s & \for s<0,\\
		f(s) & \for s \in [0,1],\\
		f(1) + f'(1)(s-1) & \for s>1,
	\end{array}
	\right.
	$$
and consider the system
	\begin{align}\nonumber
		&\io (s_t\phi(x) + k(c^P) \tilde f'(s) \nabla s\cdot\nabla \phi(x))\dd x + \ipo \alpha(x)(\tilde f(s) - \tilde f(s^{\partial\Omega}))\phi(x) \dd 
S(x)\\ \label{e1}
		&\qquad = \io Q_R(h) s Q_R(1-s)\phi(x)\dd x,\\ \nonumber
		&\io (h_t\psi(x) + (s\nabla {h} - h v^R)\cdot\nabla\psi(x))\dd x - \ipo 
\beta(x)s (h^{\partial\Omega} - h)^+\psi(x)\dd S(x)\\ \label{e2}
		&\qquad = - \io {h} s(1-s) \psi(x)\dd x
	\end{align}
for all $\phi, \psi \in V$. We first construct and solve in Section \ref{time} a time-discrete approximating system of \eqref{e1}--\eqref{e2}, and 
derive estimates independent of the time step. In Section \ref{limi}, we let the time step tend to $0$ and prove that the limit is a solution $(s,h)$ to \eqref{e1}--\eqref{e2}. We also prove that this solution has the property that $s \in [0,1]$, $h$ is positive and bounded, and $v$ is bounded, so that for $R$ sufficiently large, the truncations are never active and the solution thus satisfies \eqref{e01}--\eqref{e02} as well.

	
	\section{Time discretization}\label{time}
	
	For proving Theorem \ref{t1}, we first choose $n \in \nat$ and replace \eqref{e1}--\eqref{e2} with the following time-discrete system with time step $\tau = \frac{T}{n}$:
	\begin{align}\nonumber
		&\io \left(\frac1\tau (s_j - s_{j-1})\phi(x) + k(c^P_{j-1}) \tilde f'(s_j) \nabla s_j\cdot\nabla \phi(x)\right)\dd x + \ipo \alpha(x)(\tilde f(s_j) - \tilde f(s^{\partial\Omega}_j))\phi(x) \dd S(x) \\ \label{de1}
		&\qquad = \io Q_R(h_{j-1}) s_j Q_R(1-s_j)\phi(x)\dd x,\\ \nonumber
		&\io \left(\frac1\tau(h_j {-} h_{j-1})\psi(x) + (s_{j-1}\nabla {h_j} {-} h_j v^R_{j-1})\cdot\nabla\psi(x)\right)\dd x - \ipo \beta(x)s_{j-1} (h^{\partial\Omega}_j - h_j)^+\psi(x)\dd S(x)\\ \label{de2}
		&\qquad = - \io {h_j} s_{j-1}(1-s_{j-1}) \psi(x)\dd x,
	\end{align}
	for $j=1, \dots, n$ with initial conditions $s_0 = s^0$, $h_0 = h^0$, and with $c^P_i = 0$ for $i\le 0$, with
\begin{align}\label{qj}
q_j(x) &= -k(c^P_j(x)) \nabla \tilde f(s_j(x)) & \for\ x \in \Omega,\\ \label{vj}
v_j(x) &= \io \sigma(x-y) q_j(y) \dd y & \for\ x \in \Omega, \\ \label{sj}
s^{\partial\Omega}_j(x) &= \frac{1}{\tau}\int_{t_{j-1}}^{t_j} s^{\partial\Omega}(x,t)\dd t & \for\ x \in \partial\Omega,\\ \label{hj}
h^{\partial\Omega}_j(x) &= \frac{1}{\tau}\int_{t_{j-1}}^{t_j} h^{\partial\Omega}(x,t)\dd t & \for\ x \in \partial\Omega,
\end{align}
where $t_j = j\tau$ for $j=0,1,\dots, n$. Moreover, we define inductively
\be{cond} c^P_j - c^P_{j-1} = \tau h_j s_j (1-s_j) \ \mbox{ for } \ j=1,\dots, n.
\ee
We now prove the existence of solutions to \eqref{de1}--\eqref{de2} and derive a series of estimates which will allow us to pass to the limit as $n \to \infty$. We denote by $C$ any positive constant depending possibly on the data and independent of $n$.
	
	For $\ve >0$ we denote by $H_\ve: \real \to \real$ the function
	\be{f1}
	H_\ve(r) = \left\{
	\begin{array}{ll}
		0 &\for r\le 0,\\[2mm]
		\frac{r}{\ve} & \for r \in (0,\ve),\\[2mm]
		1 & \for r\ge \ve
	\end{array}
	\right.
	\ee
	as a Lipschitz continuous regularization of the Heaviside function, and by $\hat H_\ve$ its antiderivative
	\be{f2}
	\hat H_\ve(r) = \left\{
	\begin{array}{ll}
		0 &\for r\le 0,\\[2mm]
		\frac{r^2}{2\ve} & \for r \in (0,\ve),\\[2mm]
		r - \frac{\ve}{2} & \for r\ge \ve
	\end{array}
	\right.
	\ee
	as a continuously differentiable approximation of the ``positive part'' function. Note that we have $r H_\ve(r) \le 2 \hat H_\ve(r)$ for all $r \in \real$.

\begin{lemma}\label{l1}
		Let Hypothesis \ref{h1} hold. Then for all $n$ sufficiently large there exists a solution $h_j, \, s_j$ of the time-discrete system \eqref{de1}--\eqref{de2} with initial conditions $s_0 = s^0$, $h_0 = h^0$, $c^P_i = 0$ for $i\le 0$, which satisfies the bounds:
\be{dest1}
	s^\flat \le s_j(x) \le 1 \quad a.\,e.\quad  \for j=0,1,\dots, n.
	\ee
\be{dest3}
	h_j(x) \ge 0 \quad a.\,e.\quad \for j=0,1,\dots, n.
	\ee
	\end{lemma}

\begin{proof} 
To prove the existence, we proceed by induction. Assume that the solution to \eqref{de1}-\eqref{de2} is available for $i=1, \dots, j-1$ with the properties \eqref{dest1}--\eqref{dest3}. Then Eq.~\eqref{de1} for the unknown $s:= s_j$ is of the form
\be{mono}
\begin{aligned}
	&\io (a_0(x,s)\phi(x) + a_1(x)\tilde{f}'(s)\nabla s \cdot \nabla\phi(x))\dd x + \ipo \alpha(x)(\tilde f(s) - a_2(x))\phi(x) \dd S(x) \\
	&= \io a_3(x)\phi(x)\dd x,
\end{aligned}
\ee
where
$$
a_0(x,s) = \frac1\tau s(x) - Q_R(h_{j-1}(x)) s(x) Q_R(1-s(x))
$$
and $a_k$, $k=1,2,3$, are given functions which are known from the previous step $j-1$. For $n>TR^2$ the function $s \mapsto a_0(x,s)$ is increasing. Hence, \eqref{mono} is a monotone elliptic problem, and a unique solution exists by virtue of the Browder-Minty Theorem, see \cite[Theorem 10.49]{rr}.
Similarly, Eq.~\eqref{de2} is for the unknown function $h := h_j$ of the form
\begin{align}\nonumber
&\io (a_4(x)h\psi(x) + (a_5(x)\nabla h - a_6(x) h) \cdot \nabla\psi(x))\dd x -\ipo \beta(x) a_7(x)(a_8(x)- h)^+\psi(x)\dd S(x)\\ \label{mono2}
&\qquad = \io a_9(x)\psi(x)\dd x,
\end{align}
which we can solve in an elementary way in two steps. First, we consider the PDE
\begin{align}\nonumber
&\io (a_4(x)h\psi(x) + (a_5(x)\nabla h - a_6(x) w) \cdot \nabla\psi(x))\dd x -\ipo \beta(x) a_7(x)(a_8(x)- h)^+\psi(x)\dd S(x)\\ \label{mono3}
&\qquad = \io a_9(x)\psi(x)\dd x
\end{align}
with a given function $w \in L^2(\Omega)$. Here again the functions $a_k$, $k=4,\dots,9$, are known. We find a solution $h$ to \eqref{mono3} once more by the Browder-Minty Theorem. Since  $a_4(x) \ge \frac1\tau$, $a_5(x) = s_{j-1}(x) \ge s^\flat$, and $a_6(x) = v_{j-1}^R \in [-R,R]$, we see that for $n > TR^2/2s^\flat$, the mapping which with $w$ associates $h$ is a contraction on $L^2(\Omega)$, and the solution to \eqref{mono2} is obtained from the Banach Contraction Principle.

To derive the bounds for the solution, we first test \eqref{de1} by $\phi = H_\ve(s_j-1)$ (or any other function of the form $g(s_j - 1)$ with $g$ Lipschitz continuous, nondecreasing, and such that $g(s) = 0$ for $s \le 0$). The right-hand side identically vanishes, whereas the boundary term and $\nabla s_j\cdot \nabla H_\ve(s_j-1)$ are nonnegative, which yields that
$$
\io (s_j - s_{j-1})H_\ve(s_j - 1)\dd x \le 0.
$$
From the convexity of $\hat H_\ve$ we obtain that $(s_j - s_{j-1})H_\ve(s_j - 1) \ge \hat H_\ve (s_j-1)- \hat H_\ve (s_{j-1}-1)$, hence,
$$
\io \hat H_\ve (s_j-1) \dd x \le \io \hat H_\ve (s_{j-1}-1) \dd x.
$$
We have by hypothesis $s_{j-1} \le 1$ a.\,e., and by induction we get
\be{dest1-r}
s_j(x) \le 1 \quad a.\,e.\quad  \for j=0,1,\dots, n.
\ee
We further test \eqref{de1} by $\phi = -H_\ve(s^\flat-s_j)$. Then both the boundary term and the elliptic term give a nonnegative contribution, and using again the convexity of $\hat H_\ve$ we have
	\begin{align*}
		\frac1\tau \io (\hat H_\ve (s^\flat-s_j)  - \hat H_\ve (s^\flat-s_{j-1})) \dd x &\le
		\io -s_j H_\ve(s^\flat-s_j) Q_R(h_{j-1})Q_R(1-s_{j})  \dd x\\
		&\le \io (s^\flat-s_j) H_\ve(s^\flat-s_j) Q_R(h_{j-1}) Q_R(1-s_j) \dd x\\
		&\le 2R^2\io \hat H_\ve (s^\flat-s_j)\dd x,
	\end{align*}
	and from the induction hypothesis we get for $n >2TR^2$ that
	\be{dest2}
	s_j(x) \ge s^\flat \quad a.\,e.\quad \for j=0,1,\dots, n.
	\ee
	We have in particular $Q_R(1-s_j) = 1-s_j$ for $R\ge 1$ as well as $\tilde{f} = f$.

Test \eqref{de2} by $\psi = -H_\ve(-h_j)$. Then
	\begin{align*}
		&\frac{1}{\tau} \io (\hat H_\ve (-h_j) - \hat H_\ve (-h_{j-1}))\dd x + \io s^\flat H'_\ve(-h_j)|\nabla h_j|^2 \dd x
		\le \io h_j H'_\ve(-h_j)v^R_{j-1}\cdot \nabla h_j \dd x\\
		&\qquad \le \frac{s^\flat}2\io  H'_\ve(-h_j)|\nabla h_j|^2 \dd x + \frac1{2s^\flat} \io h_j^2 |v^R_{j-1}|^2 H'_\ve(-h_j)\dd x
	\end{align*}
	We have $0 \le h_j^2 H'_\ve(-h_j) \le \ve$ and $\lim_{\ve \to 0} \hat H_\ve(-h_j) = (-h_j)^+$, hence, passing to the limit as $\ve \to 0$, by induction we get \eqref{dest3}.
\end{proof}

\begin{lemma}\label{2}
		Let Hypothesis \ref{h1} hold. Then $h_j$ satisfies the following estimate
\be{dest4}
	\io h_j(x)\dd x \le C
	\ee
	independently of $j=0,1, \dots, n$.
\end{lemma}

\begin{proof}
Test \eqref{de2} by $\psi = 1$. Note that the boundary term is bounded above by a multiple of $h^\sharp$ and the right-hand side is negative or zero, so that
$$
\frac{1}{\tau} \io (h_j - h_{j-1}) \dd x  \le \hat{C},
$$
that is,
$$
\io h_j\dd x  \le \tau\hat{C} + \io h_{j-1}\dd x.	
$$
Summing up over $j=1, \dots, j^*$ we get
$$
\io h_{j^*}(x)\dd x \le \io h_{0}(x)\dd x + T\hat{C},
$$
for an arbitrary $j^*$, which yields \eqref{dest4}.
\end{proof}
	
The main issue will be a uniform upper bound for $h_j$ which will be 
obtained by a time-discrete variant of the Moser-Alikakos iteration technique presented in \cite{ali}. We start from some preliminary integral estimates of $s_j$.

\begin{lemma}\label{3}
Let Hypothesis \ref{h1} hold. Then $s_j$ satisfy the bounds
\be{dest5}
	\tau\sum_{j=1}^{n}\io |\nabla s_j(x)|^2\dd x +\tau\sum_{j=1}^{n}\ipo \alpha (x) s_j^2(x)\, dS(x)\le C,
	\ee
\be{dest6}
	\frac{1}{\tau} \sum_{j=1}^{n} \io |s_j - s_{j-1}|^2 \dd x + \max_{j=1, \dots, n} \io |\nabla s_j(x)|^2 \dd x \le  C \left(1+ \tau\sum_{j=0}^{n}\io h^2_j(x)\dd x \right).
	\ee
\end{lemma}

\begin{proof}
Test \eqref{de1} by $\phi = s_j$. From \eqref{dest4} we then obtain:
	\begin{align*}
	&\frac1{2\tau} \io (s_j^2 - s_{j-1}^2) \dd x + \frac{f^\flat}{2}\ipo \alpha (x) s_j^2(x)\, dS(x)+
k^\flat f^\flat\, \io |\nabla s_j(x)|^2\dd x \\
 &\le C\ipo \alpha (x)\, dS(x)+ \io h_{j-1}\dd x \le C.
	\end{align*}
Taking the sum with respect to $j=1,\dots, n$ yields
$$
\tau\sum_{j=1}^{n}\io |\nabla s_j(x)|^2\dd x + \tau\sum_{j=1}^{n}\ipo \alpha (x) s_j^2(x)\, dS(x) \le C + \io s_0^2 \dd x \le C,
$$
which is precisely \eqref{dest5}.

Let us prove \eqref{dest6} now. Test \eqref{de1} by $\phi = f(s_j)- f(s_{j-1})$. Then
	\begin{align} \nonumber
		&\frac{f^\flat}{\tau} \io |s_j - s_{j-1}|^2 \dd x + \frac12 \io \left(k(c^P_{j-1}) |\nabla f(s_j)|^2 - k(c^P_{j-2}) |\nabla f(s_{j-1})|^2\right)\dd x\\  \nonumber
		&\qquad\qquad + \frac12 \ipo \alpha(x) (f^2(s_j) - f^2(s_{j-1})) \dd S(x)\\  \nonumber
		& \qquad \le \frac12 \io (k(c^P_{j-1}) - k(c^P_{j-2})) |\nabla f(s_{j-1})|^2\dd x +
 \ipo \alpha(x) (f(s_j)f(s_j^{\partial\Omega}) - f(s_{j-1})f(s_{j-1}^{\partial\Omega})) \dd S(x) \\
\nonumber
  &\qquad\qquad   + \left( \ipo \alpha(x) f^2(s_{j-1})\dd S(x)\right)^{1/2}\, \left(  \ipo \alpha(x) |f(s_{j-1}^{\partial\Omega}) - f(s_{j}^{\partial\Omega})|^2 \dd S(x)\right)^{1/2}  \\ \label{dest6a}
  &\qquad\qquad +  \io h_{j-1}(f(s_{j})-f(s_{j-1}))\dd x.
	\end{align}
 The function $k$ is nonincreasing and the sequence $\{c^P_j\}$ is nondecreasing, so the first integral in the right-hand side of \eqref{dest6a} is negative. By Hypotheses \ref{h1}\,(i),(iii) and \eqref{sj} we further have
	$$
\frac1\tau \sum_{j=1}^{n}\ipo \alpha(x) |f(s_{j}^{\partial\Omega}) - f(s_{j-1}^{\partial\Omega})|^2 \dd S(x) \le C,
	$$
and from \eqref{dest5} it follows that
$$
\tau \sum_{j=1}^{n}\ipo \alpha(x) f^2(s_{j-1})\dd S(x) \le C,
$$
hence, by H\"older's inequality for sums,
$$
\sum_{j=1}^{n}\left( \ipo \alpha(x) f^2(s_{j-1})\dd S(x)\right)^{1/2}\, 
\left(  \ipo \alpha(x) |f(s_{j-1}^{\partial\Omega}) - f(s_{j}^{\partial\Omega})|^2 \dd S(x)\right)^{1/2} \le C.
$$
Applying Young's inequality to the last integral in the right hand side we similarly have
$$
\io h_{j-1}(f(s_{j})-f(s_{j-1}))\dd x \le C\tau \io h^2_{j-1}\dd x + \frac {f^\flat}{2\tau}\io |s_{j}-s_{j-1}|^2\dd x.
$$
Hence, taking the sum with respect to $j$ and using again \eqref{dest5}, we get \eqref{dest6}.
\end{proof}
	
\begin{remark} As a consequence of the definition of $v_j$ and $v^R_j$, we get for all $j=1, \dots, n$ that
	\be{dest7}
	\supess_{x\in \Omega}|v^R_j(x)| \le \supess_{x\in \Omega}|v_j(x)| \le C\left(1 + \left(\io|\nabla s_j(x)|^2\dd x\right)^{1/2}\right).
	\ee
\end{remark}

\begin{lemma}\label{3a}
Let Hypothesis \ref{h1} hold. Then $h_j$ satisfies the bound
\be{dest8}
	\max_{j=1, \dots, n}\io|h_j(x)|^2\dd x + \tau \sum_{j=1}^{n}\io|\nabla h_j(x)|^2 \dd x \le C.
	\ee
\end{lemma}

\begin{proof}
Test \eqref{de2} by $\psi = h_j$. Note that the boundary term is bounded above by a multiple of $(h^\sharp)^2$. Then
	\be{dstep8}
	\frac1{2\tau} \io (h_j^2 - h_{j-1}^2)\dd x + s^\flat \io |\nabla h_j(x)|^2 \dd x \le C + K_j
	\ee
	with $K_j := \io h_jv^R_{j-1}\cdot\nabla h_j\dd x$. The evaluation of this integral constitutes the most delicate part of the argument. For simplicity, we denote by $|\cdot|_r$ the norm in $L^r(\Omega)$ for $1 \le 
r \le \infty$. We first notice that by H\"older's inequality and \eqref{dest7} we have
	$$
	K_j \le C(1+ |h_j|_2|\nabla s_{j-1}|_2|\nabla h_j|_2).
	$$
Let us recall the Gagliardo-Nirenberg inequality for functions
$u \in W^{1,p}(\Omega)$ on bounded Lipschitzian domains $\Omega \subset \real^N$ in the form
\be{gn1}
|u|_q \le C\left(|u|_s + |u|_s^{1-\nu}|\nabla u|_p^{\nu}\right)
\ee
which goes back to \cite{gag,nir} and holds for every $s\le p \le q$ such 
that $\frac1q \ge \frac1p - \frac1N$, where
\be{gn2}
\nu = \frac{\frac1s - \frac1q}{\frac1N + \frac1s - \frac1p}. 
\ee
In our case we have
	\be{gn}
	|h_j|_2 \le C\left(|h_j|_1 + |h_j|_1^{1-\nu}|\nabla h_j|_2^{\nu}\right)
	\ee
	with $\nu = 3/5$. Hence, by \eqref{dest3} and \eqref{dest4},
	$$
	K_j \le C\left(1+ |\nabla s_{j-1}|_2|\nabla h_j|_2^{8/5}\right).
	$$
	Using H\"older's inequality once again we obtain
	$$
	\tau \sum_{j=1}^{n} K_j \le C\left(1+ \left(\tau \sum_{j=1}^{n}|\nabla s_{j-1}|_2^5 \right)^{1/5} \left(\tau \sum_{j=1}^{n} |\nabla h_j|_2^2\right)^{4/5}\right).
	$$
	We have
	$$
	\tau \sum_{j=1}^{n}|\nabla s_{j-1}|_2^5 \le \tau \max_{j=0, \dots, n} |\nabla s_j|_2^3 \sum_{j=0}^{n}|\nabla s_j|_2^2,
	$$
	and \eqref{dest5}--\eqref{dest6} together with \eqref{gn} yield
	$$
	\tau \sum_{j=1}^{n}|\nabla s_{j-1}|_2^5 \le C\left(1+ \left(\tau \sum_{j=1}^{n}|h_j|_2^2\right)^{3/2}\right) \le
	C\left(1+ \left(\tau \sum_{j=1}^{n}|\nabla h_j|_2^2\right)^{9/10}\right),
	$$
	so that
	$$
	\tau \sum_{j=1}^{n} K_j \le C\left(1+ \left(\tau \sum_{j=1}^{n}|\nabla h_j|_2^2\right)^{49/50}\right),
	$$
	and we conclude by summing up over $j$ in \eqref{dstep8} that \eqref{dest8} is true.
\end{proof}	

\begin{corollary}
	As an immediate consequence of \eqref{dest8} and of \eqref{dest6}, \eqref{dest7} we obtain
	\begin{align} \label{dest8b}
	\frac{1}{\tau} \sum_{j=1}^{n} \io |s_j - s_{j-1}|^2 \dd x + \max_{j=1, \dots, n} \io |\nabla s_j(x)|^2 \dd x \le  C,\\ \label{dest8a}
	\max_{j=1, \dots, n}\supess_{x \in \Omega}|v_j(x)| \le C\left(1 + \max_{j=1, \dots, n} \left(\io |\nabla s_j(x)|^2 \dd x\right)^{\frac 12} \right) \le C
	\end{align}
with a constant $C$ independent of $R$ and $n$.
\end{corollary}

\begin{corollary}
The following estimate is a direct consequence of the inequality \eqref{dest8}
\be{dest9}
	\sum_{j=1}^{n}\io |h_j - h_{j-1}|^2 \dd x \le  C.
	\ee
\end{corollary}	

\begin{proof}
To get it we test \eqref{de2} by $\psi = h_j - h_{j-1}$. On the left-hand side we keep the term 
$$
\frac1{\tau}\io |h_j - h_{j-1}|^2 \dd x
$$
and move all the other terms to the right-hand side. Thanks to \eqref{dest1} and \eqref{dest8a}, the right-hand side contains only quadratic terms in $h_j$, $h_{j-1}$, $\nabla h_j$, $\nabla h_{j-1}$. The boundary term can be estimated using the trace theorem, so that we get an inequality of the form 
	\be{dstep9}
	\frac1{\tau}\io |h_j - h_{j-1}|^2 \dd x \le  C\left (1+\io \left(|h_j^2| + |h_{j-1}|^2 + |\nabla h_j|^2 + |\nabla h_{j-1}|^2\right) \dd x \right ),
	\ee
and it suffices to apply \eqref{dest8}.
\end{proof}

The next lemma shows global boundedness of $h_j$ by means of the Moser-Alikakos iteration technique.	

\begin{lemma}
Let Hypothesis \ref{h1} hold. Then $h_j$ satisfies the bound:
\be{dest10}
	\max_{j=1, \dots, n}\supess_{x \in \Omega}|h_j(x)| \le C.
	\ee
\end{lemma}
\medskip

\begin{proof}
Consider a convex increasing function $g:[0,\infty) \to [0,\infty)$ with linear growth, and test \eqref{de2} by $\psi = g(h_j)$. We define
	$$
	G(h) = \int_0^h g(u)\dd u, \quad \Gamma(h) = \int_0^h \sqrt{g'(u)} \dd u.
	$$
	This yields
	$$
	\frac1\tau \io (G(h_j) - G(h_{j-1}))\dd x + s^\flat \io |\nabla \Gamma(h_j)|^2\dd x \le Cg(C) + C\io h_j g'(h_j)|\nabla h_j| \dd x,
	$$
	hence,
	\be{ge1}
	\max_{j=1, \dots, n} \io G(h_j) \dd x + \tau \sum_{j=1}^{n} \io |\nabla \Gamma(h_j)|^2\dd x \le  G(C) + Cg(C) +  C\tau \sum_{j=1}^{n}\io h_j^2 g'(h_j) \dd x
	\ee
	with a constant $C$ independent of $n$ and of the choice of the function 
$g$. We now make a particular choice $g = g_{M,k}$ depending on two parameters $M>1$ and $k>0$, namely
	$$
	g_{M,k}(h) = \left\{
	\begin{array}{ll}
		\frac{1}{2k+1}h^{2k+1} & \for 0\le h \le M,\\[2mm]
		\frac{1}{2k+1}M^{2k+1} + M^{2k}(h-M) & \for h>M.
	\end{array}
	\right.
	$$
	Then
	\begin{align*}
		g_{M,k}'(h) &= \min\{h, M\}^{2k} = \left\{
		\begin{array}{ll}
			h^{2k} & \for 0\le h \le M,\\[2mm]
			M^{2k} & \for h>M,
		\end{array}
		\right.\\[3mm]
		G_{M,k}(h) &= \left\{
		\begin{array}{ll}
			\frac{1}{(2k+2)(2k+1)}h^{2k+2} & \for 0\le h \le M,\\[2mm]
			\frac{1}{(2k+2)(2k+1)}M^{2k+2}+ \frac{1}{2k+1}M^{2k+1}(h-M) +\frac{1}2 
M^{2k}(h-M)^2 & \for h>M,
		\end{array}
		\right.\\[3mm]
		\Gamma_{M,k}(h) &= \left\{
		\begin{array}{ll}
			\frac{1}{k+1} h^{k+1} & \for 0\le h \le M,\\[2mm]
			\frac{1}{k+1} M^{k+1}+ M^{k}(h-M)& \for h>M.
		\end{array}
		\right.
	\end{align*}
	Note that for all $h\ge 0$, $M>0$ and $k > 0$ we have
	\be{ge1-1}
	G_{M,k}(h) \le \Gamma_{M,k}^2(h) \le 4 G_{M,k}(h), \ \ h^2 g_{M,k}'(h) \le (k+1)^2 \Gamma_{M,k}^2(h), \ \ h g_{M,k}(h) \le (k+1) \Gamma_{M,k}^2(h).
	\ee
	It thus follows from \eqref{ge1} that
	\begin{align}\nonumber
		&\max_{j=1, \dots, n} \io \Gamma_{M,k}^2(h_j) \dd x + \tau \sum_{j=1}^{n} \io |\nabla \Gamma_{M,k}(h_j)|^2\dd x\\ \label{ge2}
		&\qquad \le C \left((k+1)^2 \Gamma_{M,k}^2(C)+\tau \sum_{j=1}^{n}\io h_j^2 g_{M,k}'(h_j) \dd x\right)
	\end{align}
	with a constant $C$ independent of $k$ and $M$.
	
	We now apply again the Gagliardo-Nirenberg inequality \eqref{gn1} in the form
	$$
	|u|_q \le C\left(|u|_2 + |u|_2^{1-\nu}|\nabla u|_2^\nu \right)
	$$
	to the function $u=\Gamma_{M,k}(h_j)$, with $q = 10/3$ and $\nu = 3/5$. From \eqref{ge1-1}--\eqref{ge2} it follows that
	\begin{align}\nonumber
		&\frac{1}{k+1} \left(\tau \sum_{j=1}^{n}\left|h_j \sqrt{g_{M,k}'(h_j)}\right|_q^q\right)^{1/q} \le \left(\tau \sum_{j=1}^{n}|\Gamma_{M,k}(h_j)|_q^q\right)^{1/q}\\ \label{ge3}
		&\qquad \le C\max\left\{(k+1)^2\Gamma_{M,k}^2(C),\tau \sum_{j=1}^{n}\io h_j^2 g_{M,k}'(h_j) \dd x \right\}^{1/2}.
	\end{align}
	Let us start with $k=0$. The right-hand side of \eqref{ge3} is bounded 
independently of $M$ as a consequence of \eqref{dest8}. We can therefore let $M\to \infty$ in the left-hand side of \eqref{ge3} and obtain
	$$
	\tau \sum_{j=1}^{n}|h_j|_q^q < \infty.
	$$
	We continue by induction and put $\omega_i := (q/2)^i$ for $i \in \nat$. Assuming that
	$$
	\tau \sum_{j=1}^{n}|h_j|_{2\omega_i}^{2\omega_i} < \infty
	$$
	for some $i \in \nat$ (which we have just checked for $i=1$) we can estimate the right-hand side of \eqref{ge3} for $k= \omega_i -1$ independently of $M$, and letting $M\to \infty$ in the left-hand side we conclude that
	\begin{align*}
 \frac{1}{\omega_i}\left(\tau\sum_{j=1}^n|h_j|_{2\omega_{i+1}}^{2\omega_{i+1}}\right)^{\omega_i/2\omega_{i+1}} &\le C\max\left\lbrace\omega_i^2\Gamma_{M,\omega_i-1}^2(C),\tau\sum_{j=1}^{n}|h_j|_{2\omega_i}^{2\omega_i}\right\rbrace^{1/2}\\
		&\le C\max\left\lbrace C^{2\omega_i},\tau\sum_{j=1}^{n}|h_j|_{2\omega_i}^{2\omega_i} \right\rbrace^{1/2}.
	\end{align*}
	which implies that
\be{ge4}
\left(\tau\sum_{j=1}^{n}|h_j|_{2\omega_{i+1}}^{2\omega_{i+1}}\right)^{1/{2\omega_{i+1}}} \le (C\omega_i)^{1/\omega_i}\max\left\lbrace C,\left(\tau\sum_{j=1}^{n}|h_j|_{2\omega_i}^{2\omega_i}\right)^{1/2\omega_i} \right\rbrace
\ee
with a constant $C>0$ independent of $n$ and $i$. For $i\in \nat$ set
$$
X_i := \max\left\lbrace C,\left(\tau \sum_{j=1}^{n}|h_j|_{2\omega_i}^{2\omega_i}\right)^{1/2\omega_i}\right\rbrace,
$$
and $\Lambda_i = \log X_i$. From \eqref{ge4} it follows that $\Lambda_{i+1} \le \frac{1}{\omega_i}\log(C\omega_i) + \Lambda_i$. Summing up over $i \in \nat$ we obtain
	\be{ge5}
	\left(\tau \sum_{j=1}^{n}|h_j|_{2\omega_i}^{2\omega_i}\right)^{1/2\omega_i} \le \tilde{C}
	\ee
	with a constant $\tilde{C}$ independent of $i$. The statement now follows in a standard way.
For $\ve>0$ and $j=1, \dots, n$ put $\Omega_{j,\ve}:=  \{x \in\Omega : |h_j(x)|\geq \tilde{C}+\ve\}$, where $\tilde{C}$ is the constant from \eqref{ge5}. Then
$$
\io |h_j(x)|^{2\omega_i}\dd x \ge |\Omega_{j,\ve}|(\tilde{C}+\varepsilon)^{2\omega_i},
$$
so that
$$
\tilde{C}^{2\omega_i} \ge \tau \sum_{j=1}^{n}\io|h_j(x)|^{2\omega_i}\dd 
x \ge \left(\tau \sum_{j=1}^{n}|\Omega_{j,\ve}|\right)(\tilde{C}+\varepsilon)^{2\omega_i}.
$$
Letting $i \to \infty$ we thus obtain
$$
\tau \sum_{j=1}^{n}|\Omega_{j,\ve}| \le \lim_{i\to\infty} \left(\frac{\tilde{C}}{\tilde{C}+\varepsilon}\right)^{2\omega_i} = 0.
$$
Passing to the limit as $\varepsilon \to 0$, we obtain \eqref{dest10}.
\end{proof}

	
	\section{Limit as $n \to \infty$}\label{limi}
	
We now construct piecewise linear and piecewise constant interpolations of the sequences $\{s_j\}$, $\{h_j\}$ constructed in Section \ref{time}. Since we plan to let the discretization parameter $n$ tend to $\infty$, we denote them by $\{s_j\on\}$, $\{h_j\on\}$ to emphasize the dependence on $n$.
For $x \in \Omega$ and $t\in ((j-1)\tau, j\tau]$, $j=1, \dots, n$ set
	\be{le1}
	\begin{aligned}
		\bar s\on(x,t) &= s_j\on(x),\\[2mm]
		\underline s\on(x,t) &= s_{j-1}\on(x),\\[2mm]
		\hat s\on(x,t) &= s_{j-1}\on(x) + \frac{t- (j-1)\tau}{\tau} (s_j\on(x) - s_{j-1}\on(x)),
	\end{aligned}
	\ee
	and similarly for $\bar h\on, \hat h\on, \underline h\on, \underline v\on, \underline c^{P,(n)}, \bar s^{\partial\Omega, (n)}, \bar h^{\partial\Omega, (n)}$ etc.
	
	By virtue of the above estimates we can choose $R$ sufficiently large, so that the truncations are never active, and we can rewrite the system \eqref{de1}--\eqref{de2} in the form
	\begin{align}\nonumber
		&\io (\hat s\on_t\phi(x) + k(\underline c^{P,(n)}) f'(\bar s\on) \nabla 
\bar s\on\cdot\nabla \phi(x))\dd x + \ipo \alpha(x)(f(\bar s\on) - f(\bar 
s^{\partial\Omega, (n)}))\phi(x) \dd S(x)\\ \label{ne1} &\qquad = \io \underline h\on \bar s\on (1-\bar s\on)\phi(x)\dd x,\\ \nonumber
		&\io (\hat h\on_t\psi(x) + (\underline s\on\nabla \bar h\on - \bar h\on 
\underline v\on)\cdot\nabla\psi(x))\dd x - \ipo \beta(x)\underline s ({\overline h}^{\partial\Omega, (n)} - \bar h\on)^+\psi(x)\dd S(x)\\ \label{ne2}
		&\qquad = - \io {\bar h\on} \underline s\on(1-\underline s\on) \psi(x)\dd x.
	\end{align}
	The estimates derived in Section \ref{time} imply the following bounds independent of $n$:
	\begin{itemize}
		\item $\hat h\on$ are bounded in $L^2(0,T; V)$;
		\item $\hat s\on$ are bounded in $L^\infty(0,T; V)$;
		\item $\hat s\on_t$ are bounded in $L^2(\Omega\times (0,T))$,
		\item $\hat h\on_t$ are bounded in $L^2(0,T; V')$;
		\item $\hat h\on, \hat s\on$ are bounded in $L^\infty(\Omega\times (0,T))$.
	\end{itemize}
The bound for $\hat h\on_t$ in $L^2(0,T; V')$ follows by comparison in \eqref{ne2}. Indeed, choosing arbitrary test functions $\psi \in V$ and $\zeta \in L^2(0,T)$ in \eqref{ne2}, we obtain from the above estimates that the inequality
$$
\io \hat h\on_t(x,t)\psi(x)\zeta(t) \dd x \le C(1 + |\nabla \bar h\on(t)|_2)|\psi|_V|\zeta(t)|
$$
holds for a.\,e. $t \in (0,T)$. Integrating over $(0,T)$ and owing to estimate \eqref{dest8}, we obtain the assertion.

By the Aubin-Lions Compactness Lemma (\cite[Theorem 5.1]{li}) we can find a subsequence (still labeled by $(n)$ for simplicity) and functions $s,h$ such that
	\begin{itemize}
		\item $\hat h\on\to h$, $\hat s\on\to s$ strongly in $L^p(\Omega\times (0,T))$ for every $1\le p < \infty$;
		\item $\nabla \hat h\on\to \nabla h$ weakly in $L^2(\Omega\times (0,T); \real^3)$;
		\item $\nabla\hat s\on\to \nabla s$
		weakly-* in $L^\infty(0,T; L^2(\Omega; \real^3))$.
	\end{itemize}
In fact, the Aubin-Lions Lemma guarantees only compactness in $L^2(\Omega\times (0,T))$. Compactness in $L^p(\Omega\times (0,T))$ for $p>2$ follows from the fact that the functions are bounded in $L^\infty(\Omega\times (0,T))$ by a constant $K>0$, so that, for example,
$$
\int_0^T \io |\hat s\on(x,t) - s(x,t)|^p \dd x\dd t \le (2K)^{p-2}\int_0^T \io |\hat s\on(x,t) - s(x,t)|^2 \dd x\dd t.
$$
	Note that for $t\in ((j-1)\tau, j\tau]$ we have
	$$
	|\bar s\on(x,t) - \hat s\on(x,t)| \le |s_j\on(x) - s_{j-1}\on(x)|,
	$$
	hence,
	$$
	\int_0^T\io |\bar s\on(x,t) - \hat s\on(x,t)|^2\dd x\dd t
	\le \tau \sum_{j=1}^n \io |s_j\on(x) - s_{j-1}\on(x)|^2 \dd x \le C\tau^2
	$$
	by virtue of \eqref{dest6} and \eqref{dest8}. Similarly,
	$$
	\int_0^T\io |\bar h\on(x,t) - \hat h\on(x,t)|^2\dd x\dd t
	\le \tau \sum_{j=1}^n \io |h_j\on(x) - h_{j-1}\on(x)|^2 \dd x \le C\tau
	$$
	by virtue of \eqref{dest9}. The same estimates hold for the differences
	$\underline s\on - \hat s\on$, $\underline h\on - \hat h\on$. We conclude that
	\begin{itemize}
		\item $\bar h\on\to h$, $\bar s\on\to s$ strongly in $L^p(\Omega\times (0,T))$ for every $1\le p < \infty$;
		\item $\underline h\on\to h$, $\underline s\on\to s$ strongly in $L^p(\Omega\times (0,T))$ for every $1\le p < \infty$;
		\item $\nabla \bar h\on\to \nabla h$ weakly in $L^2(\Omega\times (0,T); \real^3)$;
		\item $\nabla\bar s\on\to \nabla s$
		weakly-* in $L^\infty(0,T; L^2(\Omega; \real^3))$.
	\end{itemize}
As a by-product of the arguments in \cite[Proof of Theorem 4.2, p.~84]{ne},	we can derive the trace embedding formula
	$$
	\ipo |u|^2 \dd S(x) \le C (|u|_2^2 + |u|_2|\nabla u|_2)
	$$
	which holds for every function $u \in V$. Consequently, we obtain strong 
convergence also in the boundary terms
	\begin{itemize}
		\item $\bar h\on\big|_{\partial\Omega}\to h\big|_{\partial\Omega}$, $\bar s\on\big|_{\partial\Omega}\to s\big|_{\partial\Omega}$ strongly in $L^2(\partial\Omega\times (0,T))$.
	\end{itemize}
	We can therefore pass to the limit as $n \to \infty$ in all terms in \eqref{ne1}--\eqref{ne2} and check that $s,h$ are solutions of \eqref{e01}--\eqref{e02} modulo the physical constants provided $R$ is chosen bigger than the constants $C$ in \eqref{dest8a} and \eqref{dest10}.


\section{Numerical test}\label{nume}

\begin{figure}[htb] \label{fig1}
\begin{center}
\includegraphics[width=16.4cm]{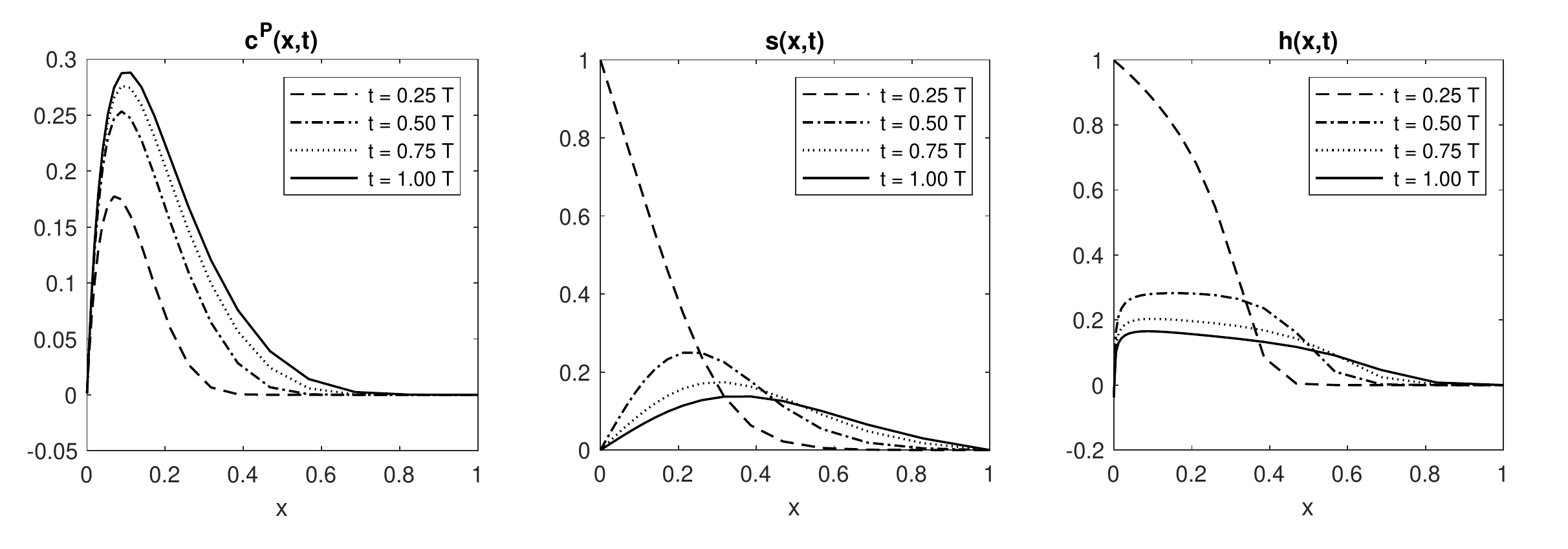}
\caption{Numerical simulations for the system \eqref{n01}--\eqref{n06}.}
\end{center}
\end{figure}

In order to illustrate the behavior of the solution, we propose a simplified 1D model with $\Omega = [0,1]$ described by the system
\begin{align} \label{n01}
\rho^W s_t(x,t) - k s_{xx}(x,t)  &= \gamma m^W {h} s(1-s),\\ \label{n02}
\rho^H h_t(x,t) - (\kappa s h_x - \rho^H hv)_x &= -\gamma m^H {h} s(1-s),
\end{align}
for $x \in (0,L)$ and $t \in (0,T)$, with boundary conditions
\begin{align} \label{n03}
k s_x(0,t) &= \alpha (s(0,t)-\bar s(0,t)),\\ \label{n04}
\kappa s h_x(0,t) - \rho^H h(0,t)v(0,t) &= -\beta s(0,t) (\bar h(0,t)- h(0,t))^+,\\ \label{n05}
k s_x(1,t) &= -\alpha s(1,t),\\ \label{n06}
\kappa s h_x(1,t) - \rho^H h(1,t)v(1,t) &= 0,
\end{align}
with some constants $\alpha>0, \beta >0$. The data are chosen so as to model the following situation: We start with initial conditions $s^0 = s^\flat$, $h^0 = 0$, and in the time interval $[0,T/4]$, we choose $\bar s = 1$ and $\bar h = 1$. This corresponds to the process of filling the structure with lime water solution until the time $t=T/4$. Then, at time $t=T/4$ we start the process of drying by switching $\bar s$ to $s^\flat$ and $\bar h$ to $0$. With these boundary data, we let the process 
run in the time interval $[T/4, T]$. Figure 1 shows the spatial distributions across the profile $x \in [0,1]$ at successive times $t=T/4, T/2, 3T/4, T$. We have chosen a finer mesh size near the origin, where the solution exhibits higher gradients. High concentration of $Ca C O_3$ near the active boundary $x=0$ exactly corresponds to the measurements shown, e.\,g., in \cite{blhvs3,sli}. The parameters of our model cannot be easily taken from the available measurements, and a complicated identification procedure would be necessary. This is beyond the scope of this paper, whose purpose is to present a model to be validated by numerical simulations. For this qualitative study we have therefore chosen fictitious parameters with simple numerical values $\rho^W = \rho^H  = m^W = m^H = \alpha = \beta = 1, s^\flat = 0, k=2\cdot 10^{-4}, \kappa=10^{-3}$ and $\gamma=10^{-2}$. The final time $T$ is determined by the number of time steps which are necessary to reach approximate equilibrium. In fact, the question of asymptotic stabilization for large times will be a subject of a subsequent study.


\section*{Acknowledgments}

The authors wish to thank Zuzana Sl\'{\i}\v zkov\'a and Milo\v s Drd\'ack\'y for stimulating discussions on technical aspects of the problem.


\end{document}